\documentclass[11pt,a4paper]{article}
\usepackage{amsfonts}
\usepackage{}

\usepackage{graphicx}
\usepackage{amsfonts,amsmath, amssymb}

\newtheorem{theo}{\textbf{\ \ \quad Theorem}}[section]

\newtheorem{lem}{\textbf{\ \ \quad Lemma}}[section]
\newtheorem{remark}{\textbf{\ \ \quad Remark}}[section]

\newtheorem{defi}{\textbf{\ \ \quad Definition}}[section]

\newcommand{\lbl}[1]{\label{#1}}

\newcommand{\be}{\begin{equation}}
\newcommand{\ee}{\end{equation}}
\newcommand\bes{\begin{eqnarray}}
\newcommand\ees{\end{eqnarray}}
\newcommand{\bess}{\begin{eqnarray*}}
\newcommand{\eess}{\end{eqnarray*}}

\newcommand{\nm}{\nonumber}

\newcommand{\ds}{\displaystyle}

\newcommand{\R}{\mathbb{R}}

{\setlength\arraycolsep{2pt}}

\title{Martingale and Weak Solutions for a Stochastic Nonlocal Burgers Equation on Bounded Intervals}

\author{Guangying Lv$^a$ and Jinqiao Duan$^b$\\
\\
\ \\
   {\small \it $^a$ Institute of Contemporary Mathematics, Henan University}\\
  {\small \it Kaifeng, Henan 475001, China}\\
  {\small \tt gylvmaths@henu.edu.cn}\\
  {\small \it $^b$ Department of Applied Mathematics, Illinois Institute of Technology}\\
  {\small \it Chicago, IL 60616 }\\
   {\small \tt duan@iit.edu }
}

\begin{document}
\maketitle

\medskip

\begin{abstract}
This work is about the existence of martingale solutions and weak solutions for a stochastic nonlocal Burgers equation on   bounded intervals. The existence of a martingale solution is shown by using a Galerkin approximation,
Prokhorov's theorem and Skorokhod's embedding theorem. The same Galerkin approximation also leads to the existence
of weak solution for the corresponding deterministic nonlocal Burgers equation on a bounded domain.



\end{abstract}

\baselineskip=15pt

\section{Introduction}
\setcounter{equation}{0}


The Fokker-Plank equation for a stochastic
differential equation with an additive Brownian motion (a Gaussian process) is a usual diffusion equation with Laplacian operator $\Delta$.
When the Brownian motion is replaced by a $\alpha$-stable L\'{e}vy motion
(a non-Gaussian process) $L_t^\alpha$, $\alpha\in(0,2)$, the Fokker-Plank
equation becomes a nonlocal partial differential equation \cite{Ap} with a nonlocal Laplacian
operator $(-\Delta)^{\frac{\alpha}{2}}$. When the drift (or vector field) of the stochastic differential equation depends on the distribution of the system evolution, this nonlocal partial differential equation becomes nonlinear.
Nonlocal Laplacian operator also appears in mathematical  models for  viscoelastic materials (e.g., Kelvin-Voigt model),  certain heat transfer processes in fractal and disordered
media, and fluid flows and acoustic propagation in porous media \cite{B,MK,NT}. Interestingly,
a nonlocal diffusion equation also arises  in pricing derivative securities in financial markets (\cite{B}).

We consider the following stochastic nonlocal Burgers equation
   \bes
   \left\{\begin{array}{llll}
du(t)=\left(-(-\Delta)^{\frac{\alpha}{2}}u-uu_x\right)dt+g(u)dW(t), \ \ t>0,x\in D,\\
u|_{D^c}=0,\\
u(x,0)=u_0(x),
  \end{array}\right.    \lbl{1.1}
  \ees
where $D=(-1,1)$ is an interval in $\mathbb{R}^1$, $D^c=\mathbb{R}^1\setminus D$,
$u_0$ is an given initial datum,  and $(-\Delta)^{\frac{\alpha}{2}}$ is the nonlocal Laplacian operator defined by the following Cauchy principal value integral
   \begin{eqnarray}  \label{Delta888}
(-\Delta)^{\frac{\alpha}{2}}u(x)=C_\alpha\int_{\mathbb{R}^1\setminus\{0\}}
\frac{u(x+y)-u(x)}{|y|^{1+\alpha}}dy, \ \ \ 0<\alpha<2,
   \end{eqnarray}
where $C_\alpha$ is a negative constant depending on $\alpha$.  The Wiener process $W_t$ will be specified later.


\medskip


\textbf{Some existing works:}
Nonlocal Burgers' type equations (deterministic or stochastic)  on the whole real line have been considered by a number of authors. For example,
Biler et al. \cite{BFW} studied the following nonlocal equation
   \bes
u_t=-(-\Delta)^{\frac{\alpha}{2}}u-uu_x, \ \ \ \  x \in \R^1,
   \lbl{1.3}\ees
and proved the existence of a unique weak solution
$u\in L^\infty(0,T;L^2(\mathbb{R}))\cap L^2(0,T;H^1(\mathbb{R}))$  for
$\alpha\in(\frac{3}{2},2]$. They further (\cite{BKW1}-\cite{BKW3}) extended this result to the equations of
the L\'{e}vy conservation laws, and obtained the
asymptotic behavior of solutions with anomalous diffusion for $1<\alpha<2$.
Bertini et al. \cite{BCJ}
studied the Burgers equation perturbed by a white noise and proved the
existence of solutions by using the Cole-Hopf transformation   in the stochastic
setting.  Wu et al. \cite{WX}, Shi and Wang \cite{ShiWang}, and  Debbi \cite{De1} considered various solutions for  a class of stochastic partial differential equations, including Burgers equation as a special case. We
remark that, in the whole space, the operator $(-\Delta)^{\frac{\alpha}{2}}$ is
similar to the Laplace operator $-\Delta$ because we can use the
Fourier transform to deal with the two operators.


\medskip

However, there are much few existing works for nonlocal Burgers' type equations (deterministic or stochastic) on bounded domains.
Mohammed-Zhang \cite{MZ} proved the existence of solutions of the stochastic Burgers
equation on a bounded interval with Dirichlet boundary conditions and anticipating initial
data by Malliavin calculus.

\medskip

\textbf{A cautious remark:}
In fact, there is another, but very different,  kind of   `fractional Laplacian operator' $(-\Delta)^{\frac{\alpha}{2}}$ on bounded domains in the literature. It is defined as a Fourier series expansion, in terms of  non-negative eigenvalues and the orthonormal basis formed by the corresponding eigenfunctions for $-\Delta$. This is similar to the textbook definition of a fractional power for a positive-definite symmetric matrix in linear algebra. For example, Debbi \cite{De} considered the fractional stochastic Navier-Stokes equations on bounded
domains with this fractional Laplacian operator.  We remark that this fractional Laplacian operator  is \emph{different} from the nonlocal Laplacian operator \eqref{Delta888} which we use here in this paper.

\bigskip

On a bounded domain the   local Laplacian operator $-\Delta$ and the nonlocal Laplacian operator $(-\Delta)^{\frac{\alpha}{2}}$ have
significant differences (see Section 2 below). Especially, the usual
fractional Sobolev spaces and embedding inequalities  will not be suitable in this context.
More information about the nonlocal operator $(-\Delta)^{\frac{\alpha}{2}}$ on bounded domains,
  are in \cite{DGLZ,CMN,GDLS,HDG}.



For the nonlocal stochastic Burgers equation (\ref{1.1}), there is no hope to use the nonlocal or anomalous diffusion  $(-\Delta)^{\frac{\alpha}{2}}$ to dominate   the convection
$uu_x$. Thus the usual method \cite{LR} and factorization method \cite{CG} are difficult to apply here.
We will adopt the method used in \cite{FG} to obtain the existence of martingale solution to
(\ref{1.1}). Because $(-\Delta)^{\frac{\alpha}{2}}$ is a nonlocal operator and
 the usual fractional Sobolev spaces \cite{NPV} will not be suitable, we will
introduce a new weighted nonlocal Sobolev space.  We also remark that the Hardy-Littlewood-Sobolev
inequality does not hold when  $\alpha$ is larger than the spatial dimension (which is $1$ in this paper).   Moreover,
we prove the existence of $L^2$ weak solution for  the nonlocal   Burgers equation (\ref{1.1}).

\medskip

The rest of this paper is organized as follows. In section 2, we will recall some results of
nonlocal Sobolev spaces.  Section 3 is concerned with the proof of the
main result on the existence of martingale solution. In section 4, we will consider the nonlocal deterministic Burgers equations in a bounded domain.

\section{Preliminaries}
\setcounter{equation}{0}

In this section, we first recall the definition of classical fractional Sobolev
space and then define a nonlocal weighted Sobolev space. Finally, we discuss some differences
between these two kinds of  spaces, and highlight
special properties for the nonlocal Sobolev spaces.

\subsection{Classical fractional Sobolev spaces }
For $s\in(0,1)$ and $p\in[1,+\infty)$, we define
   \bess
W^{s,p}(D):=\left\{u\in L^p(D):\ \frac{u(x)-u(y)}{|x-y|^{\frac{n}{p}+s}}\in L^p(D\times D)\right\};
   \eess
i.e., an intermediary Banach space between $L^p(D)$ and $W^{1,p}(D)$, endowed with
the natural norm
   \bes
\|u\|_{W^{s,p}(D)}:=\left(\int_D|u(x)|^pdx+\int_D\int_D\frac{|u(x)-u(y)|^p}{|x-y|^{n+sp}}dxdy\right)^{\frac{1}{p}},
    \lbl{2.1}\ees
where $D\subseteq\mathbb{R}^n$ is a bounded domain and the term
   \bess
[u]_{W^{s,p}(D)}:=\left(\int_D\int_D\frac{|u(x)-u(y)|^p}{|x-y|^{n+sp}}dxdy\right)^{\frac{1}{p}}
   \eess
is the so-called Gagliardo (semi) norm of $u$. Now we have the following embedding inequalities.

  \begin{lem}\lbl{l2.1}{\rm\cite[Propositions 2.1 and 2.2]{NPV}} Let $p\in[1,\infty)$ and
  $0<s\leq s'<1$. Let $D$ be an open set in $\mathbb{R}^n$ and $u:\ D\rightarrow\mathbb{R}^n$
  be a measurable function. Then
     \bess
\|u\|_{W^{s,p}(D)}\leq C\|u\|_{W^{s',p}(D)}
     \eess
for some suitable positive constant $C=C(n,s,p)\geq1$. In particular,
$W^{s',p}(D)\subseteq W^{s,p}(D)$. Furthermore, if $D$ is an open set in $\mathbb{R}^n$ of class
$C^{0,1}$ with boundary, then
     \bess
\|u\|_{W^{s,p}(D)}\leq C\|u\|_{W^{1,p}(D)}
     \eess
for some suitable positive constant $C=C(n,s,p)\geq1$. In particular,
$W^{1,p}(D)\subseteq W^{s,p}(D)$.
   \end{lem}

\begin{lem}\lbl{l2.2}{\rm\cite[Theorems 6.7 and 8.2]{NPV}}
{\rm(i)} Let $s\in(0,1)$ and $p\in[1,\infty)$ be such that $sp<n$. Let $D\subseteq\mathbb{R}^n$
be an extension domain for $W^{s,p}$. Then there exists a positive constant $C=C(n,p,s,D)$
such that, for any $u\in W^{s,p}(D)$, we have
   \bess
\|u\|_{L^q(D)}\leq C\|u\|_{W^{s,p}(D)}
    \eess
for any $q\in[p,p^*]$; i.e., the space $W^{s,p}(D)$ is continuously embedded in $L^q(D)$ for
any $q\in[p,p^*]$, $p^*=\frac{np}{n-sp}$.

If, in addition, $D$ is bounded, then the space $W^{s,p}(D)$ is continuously embedded in
$L^q(D)$ for any $q\in[1,p^*]$.

{\rm (ii)} Let $D\subseteq\mathbb{R}^n$ be an extension domain for $W^{s,p}(D)$ with no external cusps
and let $p\in[1,\infty)$, $s\in(0,1)$ be such that $sp>n$. Then, there exists $C=C(n,p,s,D)$ such
that
   \bess
\|u\|_{C^{0,\beta}(D)}\leq C\|u\|_{W^{s,p}(D)}
   \eess
for any $u\in L^p(D)$ with $\beta:=\frac{sp-n}{p}$.
 \end{lem}

For $s\in(0,1)$ and $p\in[1,\infty)$, we say that an open set $D\subseteq\mathbb{R}^n$ is an
\emph{extension domain} for $W^{s,p}$ if there exists a positive constant $C=C(n,p,s,D)$ such
that: for every function $u\in W^{s,p}(D)$ there exists $\tilde u\in W^{s,p}(\mathbb{R}^n)$ with
$\tilde u(x)=u(x)$ for any $x\in D$ and $\|\tilde u\|_{W^{s,p}(\mathbb{R}^n)}\leq C\|u\|_{W^{s,p}(D)}$.

When $s>1$ and but not an integer, we write $s=m+\sigma$, where $m$ is an integer and $\sigma\in(0,1)$.
In this case the space $W^{s,p}(D)$ consists of those equivalence classes of functions $u\in W^{m,p}(D)$
whose distributional derivatives $D^\beta u$, with $|\beta|=m$, belong to $W^{\sigma,p}(D)$, namely
   \bess
W^{s,p}(D):=\left\{ u\in W^{m,p}(D):\ D^\beta u\in W^{\sigma,p}(D) \ {\rm for\ any}\ \beta \ s.t. \ |\beta|=m\right\}.
   \eess
This is a Banach space with  norm
  \bess
\|u\|_{W^{s,p}(D)}:=\left(\|u\|^p_{W^{m,p}(D)}+\|D^\beta u\|^p_{W^{\sigma,p}(D)}\right)^{\frac{1}{p}}.
   \eess
Clearly, if $s=m$ is an integer, the space $W^{s,p}(D)$ coincides with the usual Sobolev space $W^{m,p}(D)$.
We remark that when $s>1$,   Lemmas \ref{l2.1} and \ref{l2.2} also hold.

\subsection{Nonlocal Sobolev spaces}
In this paper, we are concerned with the case with $p=2$. In order to define
the nonlocal Sobolev space, we first decompose the operator $(-\Delta)^{\frac{\alpha}{2}}$  into two components,
and then examine it as a divergence operator.
 Assume that $D\subset\mathbb{R}^n$ is an open bounded domain.

Inspired by \cite{NPV}, we rewrite the nonlocal Laplacian operator as
   \bess
(-\Delta)^{\frac{\alpha}{2}}u(x)&=&C_\alpha\int_{\mathbb{R}^n\setminus\{0\}}
\frac{u(x+y)-u(x)}{|y|^{n+\alpha}}dy\\
&=&C'_\alpha\int_{D}
\frac{u(x)-u(y)}{|x-y|^{n+\alpha}}dy+C'_\alpha \; u(x)\int_{D^c} \frac{dy}{|x-y|^{n+\alpha}}.
   \eess
where $C'_\alpha=-C_\alpha$.
Now we examine the term with
integral over $D^c$.
Denote the shortest distance to boundary by $\delta(x)=dist(x,\partial D)$ and the longest distance to boundary by $\varrho(x)=\sup_{y \in \partial D}\{dist(x,y),\ y\in \partial D\}$. Then,
we have
   \bess
B_{\varrho(x)}^c (x):=\mathbb{R}^n\setminus B_{\varrho(x)}(x)\subseteq
D^c\subseteq\mathbb{R}^n\setminus B_{\delta(x)}(x):=  B_{\delta(x)}^c (x),
   \eess
where $B_r(x)$ denotes the sphere with radius $r$ and centered at $x$. Thus,
   \bes
 \int_{D^c} \frac{dy}{|y-x|^{n+\alpha}}&\geq& \int_{B_{\varrho(x)}^c(x)} \frac{dy}{|y-x|^{n+\alpha}}\nm \\
  &\geq& \frac{C}{\varrho(x)^\alpha},
   \lbl{2.01}\ees
and
     \bes
 \int_{D^c} \frac{dy}{|y-x|^{n+\alpha}}&\leq& \int_{B_{\delta(x)}^c(x)} \frac{dy}{|y-x|^{n+\alpha}}\nm \\
  &\leq& \frac{C}{\delta(x)^\alpha},
   \lbl{2.02}\ees
 where $C$ is a positive constant. When $n=1$ and $D=(-1,1)$, we have the following exact expression
 \bess
 \int_{D^c} \frac{dy}{|y-x|^{n+\alpha}}=\frac{1}{\alpha}\left(\frac{1}{(1+x)^\alpha}+\frac{1}{(1-x)^\alpha}\right).
   \eess
Following \cite{DGLZ}, we can get another representation of the operator. We first
give a general formula. Given the vector mapping ${\mathcal {V}(x,y)}$,
$\beta(x,y):\mathbb{R}^n\times \mathbb{R}^n\rightarrow\mathbb{R}^k$ with $\beta$ antisymmetric,
i.e., $\beta(x,y)=-\beta(y,x)$, the action of the \emph{nonlocal divergence operator} $\mathcal {D}$
on $\mathcal {V}$ is defined in \cite{DGLZ2} as
   \bess
\mathcal {D}(\mathcal {V})(x):=-\int_{\mathbb{R}^n}(\mathcal {V}(x,y)+\mathcal {V}(y,x))\cdot\beta(x,y)dy, \
\ {\rm for}\ x\in\mathbb{R}^n,
   \eess
where $\mathcal {D}(\mathcal {V}):\ \mathbb{R}^n\rightarrow\mathbb{R}$.

Given the mapping $u(x):\mathbb{R}^n\rightarrow\mathbb{R}$, the adjoint operator $\mathcal {D}^*$
corresponding to $\mathcal {D}$ is the operator whose action on $u$ is given by
   \bess
\mathcal {D}^*(u)(x,y)=-(u(y)-u(x))\beta(x,y), \ \ \ {\rm for}\ x,y\in\mathbb{R}^n,
   \eess
where $\mathcal {D}^*(u):\mathbb{R}^n\times \mathbb{R}^n\rightarrow\mathbb{R}^k$.

If $\Theta (x,y)=\Theta(y,x)$ denotes a second-order tensor satisfying
$\Theta=\Theta^T$, then we have
   \bess
\mathcal {D}(\Theta\cdot\mathcal {D}^*u)(x)=-2\int_{\mathbb{R}^n}(u(y)-u(x))\beta(x,y)\cdot
(\Theta(x,y)\cdot\beta(x,y))dx,\ \ \ {\rm for}\ x\in\mathbb{R}^n,
   \eess
where $\mathcal {D}(\Theta\cdot\mathcal {D}^*u):\mathbb{R}^n\rightarrow\mathbb{R}$. If we let
$\Theta$ be the identity matrix, and $\beta$ be such that
   \bess
2|\beta(x,y)|=\frac{C'_\alpha}{|x-y|^{n+\alpha}},
   \eess
then
      \bess
\mathcal {D}(\Theta\cdot\mathcal {D}^*u)(x)=-C'_\alpha\int_{\mathbb{R}^n}\frac{u(x)-u(y)}{|x-y|^{n+\alpha}}dy,
   \eess
where $\mathcal {D}^*(u)(x,y)=(u(x)-u(y))\frac{x-y}{|x-y|^{\frac{n+\alpha}{2}+1}}$.

In particular, for $n=1$, we obtain
    \bess
(-\Delta)^{\frac{\alpha}{2}}u=C'_\alpha\int_{\mathbb{R}}\frac{u(x)-u(y)}{|x-y|^{1+\alpha}}dy
=-\mathcal {D}(\Theta\cdot\mathcal {D}^*u)(x),
   \eess
that is, the operator $(-\Delta)^{\frac{\alpha}{2}}$ is a divergence operator. For simplicity, we will
set $C'_\alpha=1$.

Direct calculations lead to
   \bes
((-\Delta)^{\frac{\alpha}{2}}u,u)_{L^2(D)}&=&-(\mathcal {D}(\mathcal {D}^*u)(x),u(x))_{L^2(D)}\nm\\
&=&-(\mathcal {D}(\mathcal {D}^*u)(x),u(x))_{L^2(\mathbb{R}^n)}\nm\\
&=&\int_{\mathbb{R}^n}\int_{\mathbb{R}^n}\frac{|u(x)-u(y)|^2}{|x-y|^{n+\alpha}}dydx\nm\\
&=&\int_{\mathbb{R}^n}\int_{\mathbb{R}^n}\frac{2u^2(x)}{|x-y|^{n+\alpha}}dydx+
\int_{\mathbb{R}^n}\int_{\mathbb{R}^n}\frac{2u(x)u(y)}{|x-y|^{n+\alpha}}dydx\nm\\
&=&2\int_{D}u^2(x)\int_{\mathbb{R}^n}\frac{1}{|x-y|^{n+\alpha}}dydx+
\int_{D}\int_{D}\frac{2u(x)u(y)}{|x-y|^{n+\alpha}}dydx\nm\\
&=&\int_{D}\int_{D}\frac{|u(x)-u(y)|^2}{|x-y|^{n+\alpha}}dydx+
\int_{D}\int_{D^c}\frac{2u^2(x)}{|x-y|^{n+\alpha}}dydx\nm\\
&=&[u]^2_{W^{\frac{\alpha}{2},2}(D)}+\int_{D}\int_{D^c}\frac{2u^2(x)}{|x-y|^{n+\alpha}}dydx,
   \lbl{2.2}\ees
where we have used the fact that $u|_{D^c}=0$. In particular, when $n=1$ and $D=(-1,1)$, we have
   \bess
((-\Delta)^{\frac{\alpha}{2}}u,u)_{L^2(D)}=[u]^2_{W^{\frac{\alpha}{2},2}(D)}
+\frac{2}{\alpha}\int_{D}u^2(x)\left(\frac{1}{(1+x)^\alpha}
+\frac{1}{(1-x)^\alpha}\right)dx.
   \eess
We remark that $W^{\frac{\alpha}{2},2}(D)$ is a Hilbert space. Actually, a scalar product is
   \bes
(u,v)_{W^{\frac{\alpha}{2},2}(D)}=\int_D u(x)v(x)dx+
\int_{D}\int_{D}\frac{(u(x)-u(y))(v(x)-v(y))}{|x-y|^{1+\alpha}}dydx.
   \lbl{2.3}\ees
 It follows from (\ref{2.2}) that the definition of fractional Sobolev space is
not suitable. The reason is that we cannot make sure that the term
   \bess
\int_{D}\int_{D^c}\frac{2u^2(x)}{|x-y|^{n+\alpha}}dydx<\infty.
   \eess
Therefore, we will introduce a weighted nonlocal Sobolev space $W^{s,p}_\rho(D)$, with
$0<s<1$, $p\geq1$,  and a `weight' function
   \bess
\rho(x)=\int_{D^c}\frac{2}{|x-y|^{n+\alpha}}dy.
    \eess
By (\ref{2.01}) and (\ref{2.02}), we see that $\rho(x)$ has strictly
positive lower bound. When $n=1$ and $D=(-1,1)$, we have
   \bess
\rho(x)=\frac{2}{\alpha}\left(\frac{1}{(1+x)^\alpha}+\frac{1}{(1-x)^\alpha}\right).
    \eess
Define
   \bess
\|u\|_{W^{s,p}_\rho(D)}:=\left(\int_D\rho(x)|u(x)|^pdx+\int_D\int_D\frac{u(x)-u(y)}{|x-y|^{n+sp}}dxdy\right)^{\frac{1}{p}}.
   \eess
It follows from \cite{KO} that $W^{s,p}_\rho(D)$ is a Banach space.
From (\ref{2.2}), we know that $\|\mathcal {D}^*u\|_{L^2(D)}=\|u\|_{W^{s,2}_\rho(D)}$. Corresponding
to (\ref{2.3}), we can define
   \bess
{}_{W^{-s,2}_\rho(D)}\langle u,v\rangle_{W^{s,2}_\rho(D)}=\int_D\rho(x) u(x)v(x)dx+
\int_{D}\int_{D}\frac{(u(x)-u(y))(v(x)-v(y))}{|x-y|^{1+2s}}dydx,
    \eess
where $W^{-s,p}_\rho(D)$ denotes the topological dual space of $W^{s,p}_\rho(D)$.
We can verify that
Lemmas \ref{l2.1} and \ref{l2.2} also hold for the weighted nonlocal Sobolev space introduced here.
Additionally, the weighted nonlocal Sobolev space introduced here is consistent with the definition
of solution to equation (\ref{1.1}) in \cite{DGLZ}. That is, for $s=\frac{\alpha}{2}$,  we have
   \bess
\|u\|_{W^{s,p}_\rho(D)}^2=\int_{\mathbb{R}^n}\int_{\mathbb{R}^n}u(-\Delta)^{\frac{\alpha}{2}}udxdy.
 \eess
If we define
   \bess
H^s(\mathbb{R}^n)=\left\{u\in L^2(\mathbb{R}^n): \int_{\mathbb{R}^n}\int_{\mathbb{R}^n}u(-\Delta)^{\frac{\alpha}{2}}udxdy<\infty\right\},
   \eess
then
   \bess
W^{s,p}_\rho(D)=\{u\in H^s(\mathbb{R}^n),\ u\equiv0\ {\rm in}\ \mathbb{R}^n\setminus D\}.
   \eess
 Thus a weighted nonlocal Sobolev space is defined.
It is known that
  \bess
\|u\|_{W^{-s,p}_\rho(D)}=
\sup_{v\in W^{s,p}_\rho(D), \|v\|_{W^{s,p}_\rho(D)}\leq1}
{}_{{W^{-s,p}_\rho(D)}}\langle u,v\rangle_{{W^{s,p}_\rho(D)}}.
   \eess
Similar to (\ref{2.2}), we have
   \bess
{}_{W^{-\frac{\alpha}{2},2}_\rho(D)}\langle(-\Delta)^{\frac{\alpha}{2}} u,v\rangle_{{W^{\frac{\alpha}{2},2}_\rho(D)}}
&=&\langle(-\Delta)^{\frac{\alpha}{2}} u,v\rangle_{L^2(D)}\nm\\
&\leq& \|\mathcal {D}^*u\|_{L^2(D)}
\|\mathcal {D}^*v\|_{L^2(D)}\nm\\
&=&\|u\|_{W_\rho^{\frac{\alpha}{2},2}}\|v\|_{W_\rho^{\frac{\alpha}{2},2}},
    \eess
which implies that
   \bes
\|(-\Delta)^{\frac{\alpha}{2}} u\|_{W^{-\frac{\alpha}{2},2}_\rho(D)}\leq C
\|u\|_{W^{\frac{\alpha}{2},2}_\rho(D)}.
   \lbl{2.4}\ees
 When $s>1$ and $s=m+\sigma$ with $\sigma\in(0,1)$, we define
   \bess
 W^{s,p}_\rho(D)=\{u\in W^{m,p}_\rho(D), \ \mathbb{D}^\gamma u\in W^{\sigma,p}_\rho(D),\ |\gamma|=m\}
    \eess
with the norm
   \bess
\|u\|_{ W^{s,p}_\rho(D)}=\|u\|_{W^{\sigma,p}_\rho(D)}+\left(\sum_{k\leq m}\int_D\rho(x)|\mathbb{D}^ku(x)|^pdx\right)^{\frac{1}{p}}.
   \eess

Before we end the this section, we present the following remark, which shows the difference
between the classical Sobolev space and nonlocal Sobolev space, that is, a difference between
  local operators and nonlocal operators.

\begin{remark}\lbl{r2.1} {\rm 1}. From the above definitions of  two kinds of Sobolev
spaces, it is clear that $W^{s,p}_\rho(D)\subset W^{s,p}(D)$.
In particular, it follows from Lemma {\rm \ref{l2.1}} that
  \bess
W^{\frac{\alpha}{2},2}_\rho(D)\subset W^{\frac{\alpha}{2},2}(D)\subset L^2(D).
  \eess

{\rm2}.   For
$u\in W^{s,p}(D)$,   we do not know any information about the
function $u$ on the boundary. Even if we consider the space $W_0^{s,p}(D)$,
where $W_0^{s,p}(D)=\{u\in W^{s,p}(D):\ u|_{\partial D}=0\}$, we only know   that the
function $u=0$ on the boundary and we do not know how the function $u$ becomes
$0$. For example, consider the following problem
  \bess\left\{
\begin{array}{llll}
-\Delta u=f(u), \ \ \ & {\rm in}\ D,\\
u=0,\ \ \ \ &{\rm on}\ \partial D.
   \end{array}\right.\eess
A working space for this   problem is $W^{1,2}_0(D)$. Since $-\Delta$
is a local operator, we do not know how the solution $u(x)$ becomes   $0$
when $x\rightarrow \partial D$. However, in problem {\rm(\ref{1.1})},
the operator $(-\Delta)^{\frac{\alpha}{2}}$ is a nonlocal operator, that is,
it is defined in the whole space. So, it has   information about how $u$ becomes
$0$ as $x\rightarrow \partial D$. In fact, from the definition of nonlocal Sobolev
space, we know that $\frac{u(x)}{\delta(x)}\rightarrow0$ as $x\rightarrow \partial D$, which dictates how $u$ becomes $0$ near boundary. It
  coincides  with the result of Theorem {\rm 1.2} in {\rm\cite{OS}} or also
{\rm\cite{OS1}}. This is a significant
difference between the fractional and nonlocal  Sobolev spaces.

$3$.When $D$ is teh whole space $\mathbb{R}^n$, that is, $D^c=\emptyset$, then {\rm(\ref{2.2})} becomes
   \bess
((-\Delta)^{\frac{\alpha}{2}}u,u)_{L^2(D)}=\int_{\mathbb{R}^n}\int_{\mathbb{R}^n}\frac{|u(x)-u(y)|^2}{|x-y|^{n+\alpha}}dydx,
    \eess
which coincides with the definition of $H^s(\mathbb{R}^n)$. Therefore,
our definition of nonlocal Sobolev space is quite natural.

\end{remark}

\section{Martingale solution for a stochastic nonlocal Burgers equation}
\setcounter{equation}{0}

In this section we consider  martingale solution for the stochastic nonlocal Burgers equation  \eqref{1.1}.

Let $W(t)$ be a Wiener process defined on a certain complete probability
space $(\Omega,\mathcal {F},P)$ and take values in the separable Hilbert space
$U$, with incremental covariance operator $Q$. Let $(\mathcal {F}_t)_{t\geq0}$
be the $\sigma$-algebras generated by $\{W(s),0\leq s\leq t\}$, then $W(t)$ is a martingale
relative to $(\mathcal {F}_t)_{t\geq0}$ and we have the following representation
of $W(t)$:
   \bess
W(t)=\sum_{i=1}^\infty\beta_i(t)e_i,
   \eess
where $\{e_i\}_{i\geq1}$ is an orthonormal set of eigenvectors of $Q$, $\beta_i(t)$
are mutually independent real Wiener processes with incremental covariance $\lambda_i>0$,
$Qe_i=\lambda_ie_i$ and $TrQ:=\sum_{i=1}^\infty\lambda_i<\infty$. For an operator
$G\in L_2(U,H)$, the space of all bounded linear operators from $U$ into $H$,
we denote by $\|G\|_2$ its Hilbert-Schmidt norm, i.e.,
    \bess
\|G\|_2^2:=Tr(GQG^*).
   \eess

Throughout this paper, we assume that $V=W^{\frac{\alpha}{2},2}_\rho(D)$,
$H=L^2(D)$ and $V_1=W^{-\delta,2}_\rho(D)$, $\delta>2+\alpha$. Then
we have
    \bess
V\subseteq H=H^*\subseteq V^*:=W^{-\frac{\alpha}{2},2}_\rho(D)\subseteq V_1.
   \eess
In addition, we make the following assumption.

(\textbf{C}) The noise intensity $g:\ H\rightarrow L_2(U,H)$ is continuous and
   \bess
\|g(u)\|_{L_2(U,H)}^2&\leq& C\|u\|_H^2+\lambda,\\
\|g(u)-g(v)\|_{L_2(U,H)}^2&\leq& C\|u-v\|_H^2,\ \ \ \forall u,v\in H
   \eess
for some positive real numbers $C$ and $\lambda$. Here and hereafter,
we assume $C$ is a positive constant and may be different from line to line.

\begin{defi}\lbl{d2.1} We say that there exists a martingale solution
of the equation {\rm(\ref{1.1})} if there exists a stochastic basis
$(\Omega,\mathcal {F},\{\mathcal {F}_t\},\mathbb{P})$, a
Wiener process $W$ on the space $U$ and progressively measurable process
$u:[0,T]\times\Omega\rightarrow H$, with $\mathbb{P}$-a.e. paths
   \bess
u(\cdot,\omega)\in L^\infty(0,T;H)\cap L^2(0,T;V)\cap C([0,T];V_1),
   \eess
such that $\mathbb{P}$-a.e., the identity
    \bes
&&(u,v)_H+\int_0^t{}_{V^*}\langle(-\Delta)^{\frac{\alpha}{2}}u,v\rangle_Vds+
\int_0^t{}_{V^*}\langle uu_x,v\rangle_Vds\nm\\
&=&(u_0,v)_H+{}_{V^*}\langle\int_0^tg(u)dW_s,v\rangle_V
  \lbl{2.5}\ees
holds   for all $t\in[0,T]$ and all $v\in V$.
\end{defi}

The main result in this section is the following theorem.

\begin{theo}\lbl{t2.1} Assume that $ \alpha \in (1, 2)$ and
$u_0$ be in $L^p(\Omega\rightarrow H;\mathcal {F}_0;\mathbb{P})$ for some $p\geq4$.
Then, under the assumption {\rm(\textbf{C})},  there exists a martingale solution for
the system {\rm(\ref{1.1})}.
\end{theo}

\bigskip

We now prove Theorem \ref{t2.1}. The main ingredients are
Galerkin approximations, Skorohod embedding theorem and the representation
Theorem.

We divide the proof into 3 Steps.

{\bf Step 1. Finite-dimensional approximation}

It follows from \cite{K} that the operator $(-\Delta)^{\frac{\alpha}{2}}$ is
positive and selfadjoint, with compact resolvent. We denote by
$0<\lambda_1<\lambda_2\leq\cdots$ the eigenvalues of $(-\Delta)^{\frac{\alpha}{2}}$,
and by $\phi_1,\phi_2,\cdots$ a complete orthonormal basis for
$H$, formed by the corresponding eigenvectors. Let
   \bess
\{e_1,e_2,\cdots\}\subset V
   \eess
be an orthonormal basis of $H$ and let $H_n:=span\{e_1,\cdots,e_n\}$ such that
$\{e_1,e_2,\cdots\}$ is dense in $V$. Let $P_n:V^*\longmapsto H_n$ be defined by
   \bess
P_ny:=\sum_{i=1}^n\langle y,e_i\rangle e_i, \ \ \ \ y\in V^*.
   \eess
Obviously, $P_n|_H$ is just the orthogonal projection onto $H_n$ in $H$ and we have
   \bess
{}_{V^*}\langle P_n(-\Delta)^{\frac{\alpha}{2}}u,v\rangle_V=\langle P_n(-\Delta)^{\frac{\alpha}{2}}u,v\rangle_H={}_{V^*}
\langle (-\Delta)^{\frac{\alpha}{2}}u,v\rangle_V, \
\ \ u\in V, \ v\in H_n,
   \eess
where ${}_{V^*}\langle\cdot,\cdot \rangle_V$ denotes the dualization between $V$ and its
dual space $V^*$. In the following section of the paper, we will omit the subscript.
Let $\{g_1,g_2,\cdots\}$ be an orthonormal basis of $U$ and
   \bess
W^n(t):=\sum_{i=1}^n\langle W(t),g_i\rangle g_i=\tilde P_nW(t),
   \eess
where $\tilde P_n$ is the orthogonal projection onto span $\{g_1,\cdots,g_n\}$ in $U$.

Then for each finite $n\in N$, we consider the following stochastic equation on $H_n$
   \bes
\left\{\begin{array}{lll}
du^n(t)=\left(-P_n(-\Delta)^{\alpha/2}u^n(t)+P_n(u^nu^n_x)\right)dt+P_ng(u^n)dW^n(t ),\ \ t\in[0,T],\\[2mm]
u^n(0)=P_nu_0=u_0^n.
  \end{array}\right.\lbl{3.1}\ees
Since the finite dimensional space stochastic differential equation (\ref{3.1}) has locally
Lipschitz and linear growth coefficient, the equation (\ref{3.1}) admits a unique strong
solution ($u^n(t)\in L^2(\Omega;C([0,T];H_n))$), see \cite{PR} for the details.

{\bf Step 2. A priori estimate}

By It\^{o} formula and noting $((-\Delta)^{\alpha/2}u,u)=\|u\|_V^2$, $(uu_x,u)=0$, we have
    \bess
\|u^n(t)\|_H^2&=&\|u^n(0)\|_H^2+2\int_0^t\langle -P_n(-\Delta)^{\alpha/2}u^n+P_n(u^nu^n_x),u^n\rangle ds\\
&&+2\int_0^t\langle u^n,P_ng(u^n)dW^n(s)\rangle+\int_0^t\|P_ng(u^n)\tilde P_n\|_2^2ds\\
&=&\|u^n(0)\|_H^2-2\int_0^t\langle (-\Delta)^{\alpha/2}u^n,u^n\rangle ds\\
&&+2\int_0^t\langle u^n,g(u^n)dW^n(s)\rangle+\int_0^t\|P_ng(u^n)\tilde P_n\|_2^2ds\\
&\leq&\|u^n(0)\|_H^2-2\int_0^t\|u^n\|_V^2 ds\\
&&+2\int_0^t\langle u^n,g(u^n)dW^n(s)\rangle+C\int_0^t\|u^n(s)\|_H^2ds+\lambda T,
   \eess
which, after taking expectations, yields that
   \bess
E\|u^n(t)\|_H^2+2 E\int_0^t\|u^n\|_V^2 ds\leq E\|u^n(0)\|_H^2+C\int_0^t\|u^n\|_H^2ds+\lambda T.
  \eess
By Gronwall's inequality, we have
       \bes
&&\sup_{0\leq t\leq T}E\|u^n(t)\|_H^2\leq c_1,\nm\\
&&E\int_0^T\|u^n\|_V^2 ds\leq c_2,
   \lbl{3.2}\ees
where $c_1,c_2$ are positive constants.

On the other hand, by It\^{o} formula and Yong's inequality, we obtain for $q\geq2$
   \bes
\|u^n(t)\|_H^q&=&\|u^n(0)\|_H^q+q(q-2)\int_0^t\|u^n\|_H^{q-4}\|(P_ng(u^n)\tilde P_n)^*u^n(s)\|_H^2ds\nm\\
&&+\frac{q}{2}\int_0^t\|u^n \|_H^{q-2}\left(2\langle -P_n(-\Delta)^{\alpha/2}u^n+P_n(u^nu^n_x),u^n\rangle+\|P_ng(u^n)\tilde P_n\|_H^2\right)ds\nm\\
&&+q\int_0^t\|u^n \|_H^{q-2}\langle u^n,P_ng(u^n)dW^n(s)\rangle\nm\\
&\leq& \|u_n(0)\|_H^q-\frac{q\theta}{2}\int_0^t\|u^n \|_H^{q-2}\|u^n\|_V^2 ds
+q(q-\frac{3}{2})\int_0^t(\lambda\|u^n \|_H^{q-2}+C\|u^n\|_H^q)ds\nm\\
&&+q\int_0^t\|u^n \|_H^{q-2}\langle u^n,P_ng(u^n)dW^n(s)\rangle.
  \lbl{3.3}\ees

It follows from Burkholder-Davis-Gundy's inequality that
   \bes
&&q\mathbb{E}\left[\sup_{0\leq t\leq T}\Big|\int_0^t\|u^n \|_H^{q-2}\langle u^n,P_ng(u^n)dW(s)\rangle\Big|\right]\nm\\
&\leq& 3q\mathbb{E}\left[
\left(\int_0^T\|u^n\|_H^{2q-2}\|g(u^n_s)\|_2^2ds\right)^{\frac{1}{2}}\right]\nm\\
&\leq& \frac{1}{3}\mathbb{E}\left[\sup_{0\leq t\leq T}\|u^n(t)\|_H^{q}\right]+C\left(\int_0^T \mathbb{E}\|g(u^n_s)\|_2^2ds\right)^{q/2}\nm\\
&\leq& \frac{1}{3}\mathbb{E}\left[\sup_{0\leq t\leq T}\|u^n(t)\|_H^q\right]+C\left(\int_0^T \mathbb{E}\|u^n_s\|_{L^2_H}^2ds\right)^{q/2}\nm\\
&\leq& \frac{1}{3}\mathbb{E}\left[\sup_{0\leq t\leq T}\|u^n(t)\|_H^q\right]+CE\int_0^T\|u^n\|_H^qds+C\lambda^qT,
   \lbl{3.4}\ees
where we have used Yong's inequality. By (\ref{3.3}), (\ref{3.4}) and using Gronwall's inequality, we obtain
   \bes
\mathbb{E}\sup_{0\leq t\leq T}\|u^n(t)\|_H^q+\mathbb{E}\int_0^T\|u^n \|_H^{q-2}\|u^n\|_V^2 ds
\leq C\left(\mathbb{E}\|u_n(0)\|_H^q+\lambda^q T\right),
   \lbl{3.5}\ees
where $C$ does not depend on $n$.

Inspired by \cite{De,FG}, we have the following lemma.
  \begin{lem}\lbl{l3.1} The sequence $\{u^n\}_{n=1,2,\cdots}$ of solutions of
  equation {\rm(\ref{3.1})} is uniformly bounded in the space
    \bess
L^2(\Omega,W^{\gamma,2}(0,T;W^{-\delta,2}_\rho(D)))\cap L^2(0,T;V),
   \eess
where $\delta<2+\alpha$ and $0<\gamma<\frac{1}{2}$.
  \end{lem}

{\bf Proof.} Inequality (\ref{3.2}) implies that $\{u^n\}_{n=1,2,\cdots}$ is
uniformly bounded in $L^2(0,T;V)$. Now we prove another part. We recall that the Besov-Slobodetski
space $W^{\gamma,p}(0,T;H)$ with $H$ being a Banach space, $\gamma\in(0,1)$ and
$p\geq1$, is the space of all $v\in L^p(0,T;H)$ such that
  \bess
\|u\|_{W^{\gamma,p}(0,T;H)}:=\left(\int_0^T\|u(t)\|_H^pdt+
\int_0^T\int_0^T\frac{\|u(x)-u(y)\|_H^p}{|t-s|^{1+\gamma p}}dtds\right)^{\frac{1}{p}},
    \eess
As $\{u^n(t)\}_{t\in[0,T]}$ is the strong solution of the finite dimensional stochastic
differential equation (\ref{3.1}), then $u^n(t)$ is the solution of the stochastic
integral equation
   \bess
u^n(t)=P_nu_0+\int_0^t(-(-\Delta)^{\alpha/2}u^n(s)+u^n(s)u^n_x(s))ds+\int_0^tP_ng(u^n)dW^n(s), \ \
a.s.
   \eess
for all $t\in[0,T]$. We denote by
   \bess
I_1(t)&=&\int_0^t(-(-\Delta)^{\alpha/2}u^n(s)+u^n(s)u^n_x(s))ds,\\
I_2(t)&=&\int_0^tP_ng(u^n)dW^n(s).
   \eess
We will prove that $I_1$ is uniformly bounded in
$L^2(\Omega,W^{\gamma,2}(0,T;W^{-\delta,2}_\rho(D)))$ and that $I_2$ is uniformly
bounded in $L^2(\Omega,W^{\gamma,2}(0,T;H(D)))$ for all $0<\gamma<\frac{1}{2}$.
Let $\phi\in W_\rho^{\delta,2}(D)$, similar to (\ref{2.4}), we get
   \bess
\Big|{}_{W_\rho^{-\delta,2}(D)}\langle u^nu^n_x,\phi\rangle_{W_\rho^{\delta,2}(D)}\Big|
&=&\Big|\langle u^n,u^n\phi_x\rangle_{L^2(D)}\Big|\\
&\leq&\|\phi_x\|_{L^\infty(D)}\|u^n\|_H^2.
   \eess
Since $2\delta>2+2\alpha$, by using Lemmas \ref{l2.1}-\ref{l2.2}, we have
   \bess
\|\phi_x\|_{L^\infty(D)}\leq C\|\phi_x\|_{H^{\frac{1}{2}+}}\leq C\|\phi\|_{W^{\delta,2}(D)}\leq
C\|\phi\|_{W_\rho^{\delta,2}(D)}.
   \eess
Therefore, we have
   \bess
\|u^nu^n_x\|_{W^{-\delta,2}_\rho(D)}\leq C\|u^n\|_H^2,
   \eess
which yields that
  \bes
\int_0^T\|I_1(t)\|_{W^{-\delta,2}_\rho(D)}^2dt&\leq&
C\int_0^T\int_0^t\left(\|(-\Delta)^{\frac{\alpha}{2}}u^n(s)\|^2_{W^{-\delta,2}_\rho(D)}
+\|u^nu^n_x\|^2_{W^{-\delta,2}_\rho(D)}\right)dsdt\nm\\
&\leq&C\int_0^T\int_0^t\left(\|u^n(s)\|^2_{H}
+\|u^n(s)\|^4_{H}\right)dsdt.
    \lbl{3.6}\ees
Here we use the following fact. Let $\phi\in W_\rho^{\delta,2}(D)$.
Since $(-\Delta)^{\frac{\alpha}{2}}$ is a divergence operator, we have
     \bes
\Big|{}_{W_\rho^{-\delta,2}(D)}\langle (-\Delta)^{\frac{\alpha}{2}}u,\phi\rangle_{W_\rho^{\delta,2}(D)}\Big|
&=&\Big|\langle u,(-\Delta)^{\frac{\alpha}{2}}\phi\rangle_{L^2(D)}\Big|\nm\\
&\leq&\|(-\Delta)^{\frac{\alpha}{2}}\phi\|_{L^2(D)}\|u\|_H.
   \lbl{3.7}\ees
Noting that $2\delta>4+2\alpha$, $\delta-\frac{1}{2}>2$ if $\alpha>1$. It follows
from Lemma \ref{l2.2} and classical Sobolev embedding theorem that
$\phi\in C_0^2(D)$. We remark that we take principle value in the definition of
the operator $(-\Delta)^{\frac{\alpha}{2}}$. It is easy to see that
   \bes
(-\Delta)^{\frac{\alpha}{2}}\phi(x)&=&\int_\mathbb{R}\frac{\phi(x+y)-\phi(x)}{|y|^{1+\alpha}}dy\nm\\
&=&\int_\mathbb{R}\frac{\phi(x+y)-\phi(x)-y\phi'(x)}{|y|^{1+\alpha}}dy\nm\\
&=&\int_D\phi''(\xi)|y|^{1-\alpha}dy<\infty,
  \lbl{3.8}\ees
where the value of $\xi$ is between $x$ and $y$. By using (\ref{3.7}) and
(\ref{3.8}), we can get (\ref{3.6}). Moreover, using H\"{o}lder inequality and
arguing as before, we obtain for $t\geq s>0$,
  \bes
\|I_1(t)-I_1(s)\|_{W_\rho^{-\delta,2}(D)}^2&=&
\|\int_s^t\left(\|(-\Delta)^{\frac{\alpha}{2}}u^n(s)+u^nu^n_xds\right)\|^2_{W^{-\delta,2}_\rho(D)}\nm\\
&\leq&C(t-s)\left(\int_s^t\|u^n(r)\|_H^2+\|u^n(r)\|_H^4dr\right).
  \lbl{3.9}\ees
Combining (\ref{3.6}), (\ref{3.9}) and (\ref{3.5}) with $q=4$, we have for $\gamma<\frac{1}{2}$
   \bess
&&\mathbb{E}\left(\int_0^T\|I_1(t)\|^2_{W_\rho^{-\delta,2}(D)}dt+\int_0^T\int_0^T
\frac{\|I_1(t)-I_1(s)\|_{W_\rho^{-\delta,2}(D)}^2}{|t-s|^{1+2\gamma}}dtds\right)^{\frac{1}{2}}\\
&\leq& C\mathbb{E}\left(\int_0^T\|u^n(r)\|_H^2+\|u^n(r)\|_H^4dr\right)^{\frac{1}{2}}\\
&\leq& C<\infty
  \eess

Now, we estimate the stochastic term $I_2$. Using the stochastic isometry, the contraction property
of $P_n$ and assumption (\textbf{C}), we have
   \bess
\int_0^T\mathbb{E}\|\int_0^tP_ng(u^n(s))dW^n(s)\|_H^2dt
&\leq& C\int_0^T\mathbb{E}\int_0^t\|g(u^n(s))\|_2^2dsdt\\
&\leq& C\int_0^T\mathbb{E}\int_0^t(\lambda+\|u^n(s)\|_H^2dsdt\\
&\leq&C<\infty.
  \eess
For $t\geq s>0$ and $\gamma<\frac{1}{2}$, the same ingredients above yield to
   \bess
\mathbb{E}\int_0^T\int_0^T
\frac{\|I_2(t)-I_2(s)\|_{H}^2}{|t-s|^{1+2\gamma}}dtds
&\leq& C\mathbb{E}\int_0^T\int_0^T
\frac{\int_s^t\|g(u^n(r))\|_{2}^2dr}{|t-s|^{1+2\gamma}}dtds\\
&\leq& C\mathbb{E}\sup_{t\in[0,T]}(1+\|u^n(t)\|_H^2)\int_0^T\int_0^T
|t-s|^{-2\gamma}dtds\\
&\leq&C<\infty.
    \eess
The proof of the Lemma is complete.  $\Box$

\begin{remark}\lbl{r3.1}
The reason we introduce the Besov-Slobodetski space is to control
the term $uu_x$. In the ordinary case, we can not get the compact
result, see the step 3.
\end{remark}

{\bf Step 3. Take weak limits}

 \begin{lem}\lbl{l3.2} {\rm\cite[Theorem 2.1]{FG}}
Let $B_0\subset B\subset B_1$ be Banach spaces, $B_0$ and $B_1$ reflexive,
with compact embedding of $B_0$ in $B$. Let $p\in(1,\infty)$ and $\alpha\in(0,1)$
be given. Let $X$ be the space
   \bess
X=L^p(0,T;B_0)\cap W^{\gamma,p}(0,T;B_1)
   \eess
endowed with the normal norm. Then the embedding of $X$ in $L^p(0,T;B)$ is compact.
\end{lem}

It follows from Lemmas \ref{l2.1} and \ref{l3.2}, we have
   \bess
W^{\gamma,2}(0,T;W^{-\delta,2}_\rho(D))\cap L^2(0,T;V)\hookrightarrow^{compact}
L^2(0,T;H).
  \eess
Therefore, we deduce that the sequence of laws $(\mathcal {L}(u^n))_n$ is tight
on $L^2(0,T;H)$. Thanks to Prokhorov's Theorem there exists a subsequence still
denoted $\{u^n\}$ for which the sequence of laws $(\mathcal {L}(u^n))_n$ converges
weakly in $L^2(0,T;H)$ to a probability measure $\mu$. By using Skorokhod's embedding
Theorem, we can construct a probability basis $(\Omega_*,\mathcal {F}_*,\mathbb{P}_*)$
and a sequence of $L^2(0,T;H)\cap C(0,T;W^{-\delta,2}_\rho(D))$-random variables $\{u^n_*\}$
and $u_*$ such that $\mathcal {L}u^n_*=\mathcal {L}u^n$, $\forall n\in N$, $\mathcal {L}(u_*)=\mu$
and $u_*^n\rightarrow u_*$ a.s. in $L^2(0,T;H)\cap C(0,T;W^{-\delta,2}_\rho(D))$. Moreover,
$u_*^n(\cdot,\omega)\in C([0,T];H_n)$. Thanks to Step 1 and the equality in law, we
obtain that the sequence $u_*^n$ converges weakly in $L^2(\Omega\times[0,T];V)$ and weakly-star
in $L^p(\Omega,L^\infty([0,T];H))$ to a limit $u_{**}$. It is easy to verify that $u_*=u_{**}$,
$dt\times dp$-a.e. and
   \bes
\mathbb{E}_*\sup_{t\in[0,T]}\|u_*(t)\|_{H}^p+\mathbb{E}_*\int_0^T\|u_*(t)\|_V^2dt\leq C<\infty.
  \lbl{3.10}\ees
We introduce the filtration
   \bess
(\mathcal {F}_n^*)_t:=\sigma\{u_*^n(s), s\leq t\},
   \eess
and construct (w.r.t. $(\mathcal {F}_n^*)_t$) the time continuous square integrable
martingale $(M_n(t),\ t\in[0,T])$ with trajectories in $C([0,T;H])$ by
   \bess
M_n(t):=u_*^n(t)-P_nu_0+\int_0^t(-\Delta)^{\frac{\alpha}{2}}u_*^n(s)ds-\int_0^tu_*^n(s)(u_*^n(s))_xds.
   \eess
The equality in law yields to the fact that the quadratic variation is given by
   \bess
\ll M_n\gg_t=\int_0^tP_ng(u_*^n(s))Qg(u_*^n(s))^*ds,
  \eess
where $g(u_*^n(s))^*$ is the adjoint of $g(u_*^n(s))$. We will prove $M_n(t)$ converges
weakly in $W^{-\delta,2}_\rho(D)$ to the martingale $M(t)$ for all $t\in[0,T]$, where $M(t)$
is given by
  \bess
M(t):=u_*(t)-u_0+\int_0^t(-\Delta)^{\frac{\alpha}{2}}u_*(s)ds-\int_0^tu_*(s)(u_*(s))_xds.
   \eess
It follows from Lemma \ref{l3.1} that $\|M_n(t)\|_{W^{-\delta,2}_\rho(D)}\leq C<\infty$,
where $C$ does not depend on $n$. Since $W^{-\delta,2}_\rho(D)$ is a Hilbert space, we have
$M_n(t)\rightharpoonup M(t)$, as $n\rightarrow\infty$, $t\in[0,T]$. Now we apply to the representation
Theorem \cite[Theorem 8.2]{PZ}, we infer that there exists a probability basis
$(\Omega^*,\mathcal {F}^*,\mathbb{R}^*,W^*)$ such that
   \bess
M(t)=\int_0^tg(u_*)dW^*(s).
   \eess
By using Burkhody-Davis-Gundy inequality and (\ref{3.10}), we have
   \bess
\mathbb{E}\sup_{t\in[0,T]}\|\int_0^tg(u_*(s))dW^*(s)\|_H^2
&\leq&C\mathbb{E}\int_0^T\|g(u_*(s))\|_2^2ds\\
&\leq&C(1+\mathbb{E}\sup_{t\in[0,T]}\|u_*(s)\|_H^2)<\infty.
  \eess
Furthermore, by using Lemmas \ref{l2.1}-\ref{l2.2}, (\ref{2.4}) and
$1<\alpha<2$, we get
   \bess
&&\mathbb{E}\int_0^T(\|(-\Delta)^{\frac{\alpha}{2}}u_*(s)\|_{W_\rho^{-\frac{\alpha}{2},2}(D)}
+\|(u_*^2(s))_x\|_{W_\rho^{-\frac{\alpha}{2},2}(D)})ds\\
&\leq&
C\mathbb{E}\int_0^T(\|u_*(s)\|_V
+\|u_*^2(s)\|_{W_\rho^{1-\frac{\alpha}{2},2}(D)})ds\\
&\leq& C\mathbb{E}\int_0^T(1+\|u_*(s)\|_V)ds<\infty,
   \eess
where we used the facts $u_*(s)\in C(\bar D)$ because of $\alpha>1$, and
   \bess
\|(u_*^2(s))_x\|_{W_\rho^{-\frac{\alpha}{2},2}(D)}\leq C\|u_*^2(s)\|_{W_\rho^{1-\frac{\alpha}{2},2}(D)}.
   \eess
By using Fourier transform, one can prove that the above inequality holds for $D=\mathbb{R}$.
Noting that $u_*\equiv0$ in $\mathbb{R}\setminus D$, we have the above inequality holds. Actually,
even if $u_*$ does not define in $\mathbb{R}\setminus D$, we can also get the desire result under the
condition that $D$ is \emph{an extension domain}.
Because the domain $D$ is \emph{an extension domain}, we can extend $u$ to $\mathbb{R}$ by
letting $u=0$ in $\mathbb{R}\setminus D$. We denote it by $\tilde u$, and obtain
that
   \bess
\|u\|_{W_\rho^{-\frac{\alpha}{2},2}(D)}\leq C\|\tilde u\|_{W_\rho^{1-\frac{\alpha}{2},2}(\mathbb{R})}
\leq C\|u\|_{W_\rho^{1-\frac{\alpha}{2},2}(D)}.
   \eess
Using the densely embedding $W_\rho^{\delta,2}(D)\hookrightarrow V$, we conclude
that (\ref{2.5}) holds in the $W_\rho^{\frac{\alpha}{2},2}(D)$-$W_\rho^{-\frac{\alpha}{2},2}(D)$-duality.
This complete the proof of Theorem \ref{t2.1}. $\Box$

\begin{remark}\lbl{r3.1} We can use the same method to deal with the following problem
   \bes\left\{\begin{array}{llll}
du_t=\left(-(-\Delta)^{\frac{\alpha}{2}}u-uu_x\right)dt+g(u)dW_t, \ \ t>0,x\in D,\\
u(x,t)|_{D^c}=g(x,t),\\
u(x,0)=u_0(x),
  \end{array}\right.\lbl{3.11}\ees
where $D=(0,1)$. By letting $v=u-g$, we can obtain the existence of martingale
solution of {\rm(\ref{3.11})} under the suitable assumption on $g$.
\end{remark}

\section{Weak solution for a deterministic nonlocal Burgers equation}
\setcounter{equation}{0}

In this section, we will consider the corresponding deterministic version of the equation \eqref{1.1} considered in the previous section.
That is, we consider the following deterministic nonlocal Burgers
equation on a bounded interval
  \bes\left\{\begin{array}{lllll}
u(t)+(-\Delta)^{\frac{\alpha}{2}}u + uu_x=0,\ \ \ t>0,\ x\in D,\\
u|_{D^c}=0,\\
u(x,0)=u_0(x),
   \end{array}\right.\lbl{4.1}\ees
where $D=(-1,1)$ and $D^c=\mathbb{R}^1 \setminus D$.

In section 3, we obtained the existence of martingale solution
to (\ref{1.1}). Unfortunately, we can not get the uniqueness of the martingale solution.
We are unaware of an existence of weak solution result about the nonlocal Burgers equation (\ref{4.1}) on a bounded
domain.  On the  whole space  $D=\mathbb{R}^1$, there are a lot of results for (\ref{4.1}); see, for example,
\cite{BFW,BKW1,BKW2,BKW3}.

In the following, we will adopt the similar method to section 3 to prove the existence of
$L^2$-solution of (\ref{4.1}). Firstly, we give the definition of $L^2$-solution. We will
adopt the same symbol as in section 2. Let
   \bess
C^n_\rho(D):=\left\{u\in C^n(D),\ \ \rho(x)u^{(n)}(x)\in L^\infty(D)\right\}.
   \eess

\begin{defi}\lbl{d4.1} We say that  $u$ is a  weak solution   of {\rm(\ref{4.1})} if
$u\in L^\infty(0,T;L^2(D))$
for each $T>0$, such that $u$ satisfies
   \bess
(u,\phi)+\int_0^t(u,(-\Delta)^{\alpha/2}\phi)ds-\frac{1}{2}\int_0^t(u^2,\phi_x)ds=(u_0,\phi),
   \eess
   for  each $\phi\in C^2_\rho(D)$.
   \end{defi}

\begin{theo}\lbl{t4.1}
For $u_0\in L^2(D)$,  there exists a  weak solution to equation {\rm(\ref{4.1})}.
  \end{theo}

{\bf Proof. } We will use a Galerkin approximation and Lemma \ref{l3.2} to prove
this Theorem. Similar to the proof of Theorem \ref{t2.1}, let
   \bess
\{e_1,e_2,\cdots\}\subset V
   \eess
be an orthonormal basis of $H$ and let $H_n:=span\{e_1,\cdots,e_n\}$ such that
$\{e_1,e_2,\cdots\}$ is dense in $V$. Let $P_n:V^*\longmapsto H_n$ be defined by
   \bess
P_ny:=\sum_{i=1}^n\langle y,e_i\rangle e_i, \ \ \ \ y\in V^*.
   \eess
Obviously, $P_n|_H$ is just the orthogonal projection onto $H_n$ in $H$ and we have
   \bess
{}_{V^*}\langle P_n(-\Delta)^{\frac{\alpha}{2}}u,v\rangle_V=\langle P_n(-\Delta)^{\frac{\alpha}{2}}u,v\rangle_H={}_{V^*}
\langle (-\Delta)^{\frac{\alpha}{2}}u,v\rangle_V, \
\ \ u\in V, \ v\in H_n,
   \eess
where ${}_{V^*}\langle\cdot,\cdot \rangle_V$ denotes the dualization between $V$ and its
dual space $V^*$. Then for each finite $n\in N$, we consider the following stochastic equation on $H_n$
   \bes
\left\{\begin{array}{lll}
\ds\frac{du^n(t)}{dt}=-P_n(-\Delta)^{\alpha/2}u^n(t)-P_n(u^nu^n_x),\ \ t\in[0,T],\\[2mm]
u^n(0)=P_nu_0=u_0^n.
  \end{array}\right.\lbl{4.2}\ees
Since the finite dimensional space stochastic differential equation (\ref{4.2}) has locally
Lipschitz and linear growth coefficient, the equation (\ref{3.1}) admits a unique strong
solution ($u^n(t)\in L^2(\Omega;C([0,T];H_n))$).

Multiplying (\ref{4.2}) by $u_n$ and integrating over $D\times[0,t)$, we have
  \bes
\|u^n(t)\|_H^2&=&\|u^n_0\|_H^2-\int_0^t((-\Delta)^{\frac{\alpha}{2}}u^n(s),u^n(s))_{H}ds\nm\\
&=&\|u^n_0\|_H^2-\int_0^t\|u^n(s)\|_{V}^2ds,
  \lbl{4.3}\ees
where we have used the facts $(u^nu^n_x,u^n)=0$ and (\ref{2.2}). Equality (\ref{4.3}) implies that
$\|u^n\|_H^2\leq \|u_0\|_H^2$, which yields that there exists a subsequence of $\{u^n\}$, still
denoted $\{u^n\}$, such that $u^n\rightharpoonup u$ as $n\rightarrow\infty$, where $u\in H$.

Next we prove that the sequence $\{u^n\}$ is uniformly bounded in the space
$W^{\gamma,2}(0,T;W^{-\delta,2}_\rho(D))\cap L^2(0,T;V)$. Inequality (\ref{4.3}) implies that $\{u^n\}_{n=1,2,\cdots}$ is
uniformly bounded in $L^2(0,T;V)$. As $\{u^n(t)\}_{t\in[0,T]}$
is the strong solution of the finite dimensional stochastic
differential equation (\ref{3.1}), then $u^n(t)$ is the solution of the stochastic
integral equation
   \bess
u^n(t)=P_nu_0-\int_0^t((-\Delta)^{\alpha/2}u^n(s)+u^n(s)u^n_x(s))ds
   \eess
for all $t\in[0,T]$. We denote by
   \bess
I(t)&=&-\int_0^t((-\Delta)^{\alpha/2}u^n(s)+u^n(s)u^n_x(s))ds.
   \eess
It is remarked that $\sup_{t\in[0,T]}\|u^n(t)\|_H^2\leq C$ implies that
$\sup_{t\in[0,T]}\|u^n(t)\|_H^4\leq C$. Similar to the proof of Lemma \ref{l3.1}, one can
prove that
   \bess
\int_0^T\|I(t)\|^2_{W_\rho^{-\delta,2}(D)}dt+\int_0^T\int_0^T
\frac{\|I(t)-I(s)\|_{W_\rho^{-\delta,2}(D)}^2}{|t-s|^{1+2\gamma}}dtds\leq C.
  \eess
It follows from Lemmas \ref{l2.1} and \ref{l3.2}, we have
   \bess
W^{\gamma,2}(0,T;W^{-\delta,2}_\rho(D))\cap L^2(0,T;V)\hookrightarrow^{compact}
L^2(0,T;H).
  \eess
Therefore, we deduce that the sequence $\{u^n\}$ converges to some $u^*$ in $L^2(0,T;H)$.
Due to the uniqueness of the limit, we obtain that $u=u^*$. Let $\phi\in C^2_\rho(D)$, then
we have
   \bess
(u^n,\phi)+\int_0^t(u^n,(-\Delta)^{\alpha/2}\phi)ds-\frac{1}{2}\int_0^t((u^n)^2,\phi_x)ds=(u^n_0,\phi).
   \eess
Let $n\rightarrow \infty$, we have
   \bess
(u,\phi)+\int_0^t(u,(-\Delta)^{\alpha/2}\phi)ds-\frac{1}{2}\int_0^t(u^2,\phi_x)ds=(u_0,\phi).
   \eess
This completes the proof.  $\Box$

\medskip

\noindent {\bf Acknowledgment} The first author was supported in part
by NSFC of China grants 11301146, 11171064.

 \end{document}